\documentclass[reqno,11pt]{amsproc}
\usepackage{amssymb}
\usepackage{amscd}
\usepackage{amsmath}
\usepackage{showkeys}

\theoremstyle{plain}

\newtheorem{proposition}{Proposition}

\newtheorem{theorem}{Theorem}
\numberwithin{equation}{section}

\def\Nset{\mbox{I\kern-.21em N}}
\def\RE{{\mbox{\rm I\kern-.21em R}}}
\def\ZZ{{\mbox{\sf Z\kern-.45em Z}}}
\def\vv{\kern.344em{\rule[.18ex]{.075em}{1.32ex}}\kern-.344em}

 \def\f{\varphi}

\def\<{\langle} \def\>{\rangle}

\begin{document}

\title{Relationship between
different types of inverse data for the one-dimensional
Schr\"odinger operator on a half-line.}
\author{ A. S. Mikhaylov}
\address{St. Petersburg   Department   of   V.A. Steklov    Institute   of   Mathematics
of   the   Russian   Academy   of   Sciences, 7, Fontanka, 191023
St. Petersburg, Russia and Saint Petersburg State University,
St.Petersburg State University, 7/9 Universitetskaya nab., St.
Petersburg, 199034 Russia.} \email{a.mikhaylov@spbu.ru}
\author{ V. S. Mikhaylov}
\address{St.Petersburg   Department   of   V.A.Steklov    Institute   of   Mathematics
of   the   Russian   Academy   of   Sciences, 7, Fontanka, 191023
St. Petersburg, Russia and Saint Petersburg State University,
St.Petersburg State University, 7/9 Universitetskaya nab., St.
Petersburg, 199034 Russia.} \email{v.mikhaylov@spbu.ru}


\maketitle

\begin{center}
{\bf Abstract.} We consider inverse dynamical, spectral, quantum
and acoustical scattering problems for the Schr\"odinger operator
on the half line. The goal of the paper is to establish the
connections between different types of inverse data for these
problems. The central objects which serve as a source for all
formulaes are kernels of so-called connecting operators and Krein
equations.
\end{center}

\section{Introduction}

This paper is of methodological character, its primary goal is to
show the connection of the different types of inverse data for the
Shr\"odinger operator on the half line. The idea of using the
connection of the inverse data in solving the inverse problems is
not knew, to mention \cite{B07,KKLM,AM,AMR,MM}. In our approach we
exploit the central objects of the Boundary Control method -- the
connecting operators and corresponding Krein equations, and show
that using just these two well-known objects leads to interesting
results.

The central object we will be dealing with is a wave equation on
the half-line with the potential $q\in L_{loc}^{1}\left(
\mathbb{R}_{+}\right)$:
\begin{equation}
\label{wave_eqnGL} \left\{
\begin{array}l
u_{tt}(x,t)-u_{xx}(x,t)+q(x)u(x,t)=0, \quad x>0,\ t>0,\\
u(x,0)=u_t(x,0)=0,\ u(0,t)=f(t).
\end{array}
\right.
\end{equation}
Here $f$ is an arbitrary $L^2_{loc}\left( \mathbb{R}_{+}\right) $
function referred to as a \emph{boundary control}, by $u^f$ we
denote the solution to (\ref{wave_eqnGL}). Let $T>0$ be fixed. The
dynamical inverse data is given by the \emph{response operator}
(the dynamical Dirichlet-to-Neumann map)
$\left(R^Tf\right):=u_{x}^f(0,t)$, and the inverse problem
associated with (\ref{wave_eqnGL}) is to recover $q(x),$ $0<x<T$,
by given $R^{2T}$. One of the efficient methods of solving this
problem is the Boundary Control method \cite{B07,AM,BM01}. The
\emph{control operator} and \emph{connecting operator } are
introduced by $W^Tf:=u^f(\cdot,T)$, $C^T:=\left(W^T\right)^*W^T$.
The fact that $C^T$ is expressed in terms of the inverse data
\cite{B07} plays an important role in BC method.

We also consider the spectral, and quantum and acoustical
scattering problems for the Schr\"odinger operator with the same
potential $q$ on the half-line $H=-\partial^2_x+q$ on
$L_2(0,\infty)$ with Dirichlet boundary condition $\phi(0)=0$. For
each problem we define corresponding data: spectral measure and
Titchmarsh-Weil function for the spectral problem, discrete
spectrum with norming coefficients and scattering matrix for the
scattering problem (we need to assume that the potential satisfy
some additional condition on growth at infinity); acoustical
response operator and acoustical response function (for this
problem we assume that potential is infinitely smooth and
compactly supported). Our aim will be to show the connection of
the dynamical data, which is the kernel of the response operator
$R^T$ with spectral and scattering data and connection od the
acoustical response with the scattering data. Some of the results
have been obtained in \cite{AMR,AM}, we list them for the sake of
completeness or give a different proof. The main objects which
play the key role in our considerations is the kernels of the
connecting operators and Krein equations. The central role of the
connecting operators in different inverse problems have been
pointed out in \cite{B07,B03,AMM}, in \cite{Dem} the author
studied the singular values of connecting operator for the
observation problem.

In the second section we set up the forward and inverse problems:
dynamical, spectral, quantum and acoustical scattering, and for
each of them introduce the corresponding inverse data. In the
third section we study in details the integral kernel of the
connecting operator and reveal the links with the spectral
function of Levitan \cite{L3}. In the last section we derive the
spectral and scattering representation of the response function
and explain the connection of the response function with the
Titchmarsh-Weil function see also \cite{AMR,AM}; we also derive
the scattering representation for the acoustical scattering
response function and establish its connection with the scattering
matrix.

\section{Inverse data}

\subsection{Dynamical inverse data}

For the potential $q\in L_{loc}^{1}\left( \mathbb{R}_{+}\right)$
we consider the initial boundary value problem for the 1d wave
equation with the potential (\ref{wave_eqnGL}) with $f$ be an
arbitrary $L^2_{loc}\left( \mathbb{R}_{+}\right) $ function
referred to as a \emph{boundary control}. It is known \cite{AM}
that the solution $u^f(x,t)$ of the problem (\ref{wave_eqnGL}) can
be written in terms of the integral kernel $w(x,s)$:
\begin{equation}
\label{wave_eqn_sol} u^f(x,t)=\left\{\begin{array}l
f(t-x)+\int_x^tw(x,s)f(t-s)\,ds, \quad x \leq t,\\
0, \quad x > t.\end{array}\right .
\end{equation}
where $w(x,s)$ is the unique solution to certain Goursat problem.
Fix $T>0$ and introduce the outer space of the system
(\ref{wave_eqnGL}), the space of controls:
$\mathcal{F}^T:=L_2(0,T)$.  The dynamical inverse data is given by
the \emph{response operator} (the dynamical Dirichlet-to-Neumann
map) $R^T:\mathcal{F}^T\mapsto \mathcal{F}^T$ with the domain $\{
f \in C^2([0,T]):\; f(0)=f'(0)=0\}$, acting by the rule:
\begin{equation*}
(R^Tf)(t)=u_{x}^f(0,t), \ t \in (0,T),
\end{equation*}
According to (\ref{wave_eqn_sol}) it has a representation
\begin{equation}
\label{react_rep} (R^Tf)(t)=-f'(t)+\int_0^tr(s)f(t-s)\,ds,
\end{equation}
where $r(t):=w_x(0,t)$ is called the \emph{response function}. The
natural set up of a inverse problem \cite{AM,B07,BM01} is to
recover the potential $q(x),$ $x\in (0,T)$ from $R^{2T}$, or what
is equivalent, from $r(t),$ $t\in (0,2T).$

We introduce the inner space of system (\ref{wave_eqnGL}), the
space of states: $\mathcal{H}^T:=L_2(0,T)$, so for all $0\leqslant
t\leqslant T,$ $u^f(\cdot,t)\in \mathcal{H}^T$. The \emph{control
operator} ${W}^T$ is defined by
$$
{ W}^T:\mathcal{F}^T \mapsto \mathcal{H}^T, \ { W}^T f=u^f(
\cdot,T),
$$
is bounded. From (\ref{wave_eqn_sol}) it follows that
\begin{equation*}
(W^Tf)(x)=f(T-x)+ \int_{x}^{T} w(x,\tau)f(T-\tau)\,d\tau.
\end{equation*}
It is not hard to show that $W^T$ is in fact boundedly invertible.
The \emph{connecting operator} $C^T: \mathcal{F}^T \mapsto
\mathcal{F}^T,$ plays a central role in the BC method. It connects
the outer space of the dynamical system (\ref{wave_eqnGL}) with
the inner space, and is defined by its bilinear product:
\begin{equation} \label{Ct} \left(C^Tf,g\right)_{\mathcal{F}^T}=
\left(u^f(\cdot,T),u^g(\cdot,T)\right)_{\mathcal{H}^T},\quad
C^T=(W^T)^*W^T.
\end{equation}
The invertibility of $W^T$ implies that $C^T$ is positive
definite, bounded and boundedly invertible in $\mathcal{F}^T.$ The
fact that $C^T$ is expressed in terms of the response operator is
widely used in BC-method. In \cite{AM} we have shown this for the
case of nonsmooth potential:
\begin{proposition} \label{pr3} For $q \in L_{loc}^1(0,T)$ and
$T >0$, operator $C^T$ admits the representation
\begin{equation}
\label{r-c2} ({C}^T  f)(t)=f(t)+\int_0^T c^T(t,s)f(s)\,ds\,, \
0<t<T\,,
\end{equation}
where
\begin{equation}
\label{c_t} c^T(t,s)=[p(2T-t-s)-p(t-s)],\quad p(t)=\frac{1}{2}
\int_0^t r(s)\,ds\end{equation}
\end{proposition}
We fix a function $y$ to be solution to the following Cauchy
problem:
\begin{equation*}
\left\{
\begin{array}l
 -y''+qy=\lambda y,\qquad x>0\\
y(0)=0,\quad y'(0)=1.
\end{array} \right.
\end{equation*}
Set up the \emph{special control problem}: to find a  control
$f^T$ that
\begin{equation*}
(W^Tf^T)(x)=\left\{\begin{array}l y(x),\ 0<x<T,\\
0,\ x>T.\end{array}\right.
\end{equation*}
\begin{theorem}
\label{Wv_control_teor} The control $f^T=W^{-1}y,$ which solves
the special control problem, is the solution of the Krein equation
$Cf=\frac{\sin{\sqrt{\lambda}(T-t)}}{\sqrt{\lambda}}$ in
$\mathcal{F}$, i.e., satisfies
\begin{equation*}
f^T(\tau)+\int_0^T
c^T(t,s)f^T(s)\,ds=\frac{\sin{\sqrt{\lambda}(T-t)}}{\sqrt{\lambda}},\qquad
0\leqslant t\leqslant T.
\end{equation*}
\end{theorem}
Notice that $\left(W^T\right)^*$ is a transformation operator: it
maps the solution of the perturbed problem to the solution of the
unperturbed (modulo shift by $T$):
\begin{equation}
\label{WT_wv_tr}
\left(W^T\right)^*y(\cdot,\lambda)=\frac{\sin{\sqrt{\lambda}(T-t)}}{\sqrt{\lambda}}
\end{equation}

\subsection{Spectral inverse data.}

we consider the Schr\"odinger operator
\begin{equation}
\label{Schrodinger} H=-\partial _{x}^{2}+q\left( x\right)
\end{equation}
on $L^{2}\left( \mathbb{R}_{+}\right) \,,\mathbb{R}_{+}:=[0,\infty
),$ with a real-valued locally integrable potential $q$ and
Dirichlet boundary condition at $x=0.$ For $z\in \mathbb{C}$
consider the solution
\begin{equation}
\label{Schr} \left\{
\begin{array}l
-\varphi''(x)+q(x)\varphi(x)=\lambda\varphi(x),\\
\varphi(0,z)=0,\,\, \varphi(0,z)=1.
\end{array}
\right.
\end{equation}
It is known \cite{LeSa} that there exist a spectral measure
$d\rho(\lambda)$, such that for all $f,g\in L^2(\mathbb{R}_+)$ the
Parseval identity holds:
\begin{equation}
\int_0^\infty f(x)g(x)\,dx=\int_{-\infty}^\infty
(Ff)(\lambda)(Fg)(\lambda)\,d\rho(\lambda), \label{Fourier_int_1}\\
\end{equation}
where $F: L_2(\mathbb{R}_+)\mapsto L_{2,\,\rho}(\mathbb{R})$ is a
Fourier transformation:
\begin{eqnarray}
(Ff)(\lambda)=\int_0^\infty
f(x)\varphi(x,\lambda)\,dx\label{Fourier_int_2}\\
f(x)=\int_{-\infty}^\infty
(Ff)(\lambda)\varphi(x,\lambda)\,d\rho(\lambda)\notag.
\end{eqnarray}
The so-called transformation operator transforms the solutions of
(\ref{Schr}) to the solution of (\ref{Schr}) with zero potential:
\begin{equation}
\label{Inv_trans_wave}
\left(I_s+L_s\right)\varphi(\cdot,\lambda)=\f(s,\lambda)+\int_0^s
w(x,s)\f(x,\lambda)\,dx=\frac{\sin{\sqrt{\lambda}s}}{\sqrt{\lambda}}.
\end{equation}

We assume that $\left( \ref{Schrodinger}\right) $ is limit point
case at $\infty $, that is, for each $z\in \mathbb{C}_{+}:=\{z\in
\mathbb{C}:\operatorname{Im}z>0\}$ the equation
\begin{equation}
-u^{\prime \prime }+q\left( x\right) u=zu  \label{Sch_eqnn}
\end{equation}
has a unique, up to a multiplicative constant, solution in $%
L_{2}$ at $\infty ~$, we denote this solution by $u_+$:
\begin{equation*}
\int_{\mathbb{R}_{+}}\left\vert u_{+}\left( x,z\right) \right\vert
^{2}dx<\infty ,\quad z\in \mathbb{C}_{+}.
\end{equation*}
Then the \emph{Titchmarsh-Weyl }$m$\emph{-function}, $m(z),$ is
defined for $z\in \mathbb{C}_{+}$ as
\begin{equation*}
m\left( z\right) :=\frac{u_{+}^{\prime }\left( 0,z\right)
}{u_{+}\left( 0,z\right) }.
\end{equation*}%
The function $m(z)$ is analytic in $\mathbb{C}_{+}$ and satisfies
the Herglotz property: $m:\mathbb{C}_{+}\rightarrow
\mathbb{C}_{+}$, so $m$ satisfies a Herglotz representation
theorem,
\begin{equation*}
m\left( z\right) =c+\int_{\mathbb{R}}\left( \frac{1}{t-z}-\frac{t}{1+t^{2}}%
\right) d\rho \left( t\right) ,
\end{equation*}%
where $c=\operatorname{Re}m\left( i\right) ~$and $\rho $ is
spectral measure of $H$. The measure can be recovered from $m(z)$
by the rule:
\begin{equation*}
d\rho \left( t\right) =\text{w-}\lim_{\varepsilon \rightarrow
+0}\frac{1}{\pi }\operatorname{Im}m\left( t+i\varepsilon \right)
dt.
\end{equation*}
On the problems of uniqueness and recovering the potential from
the Weyl function we refer to to classical papers by Borg
\cite{Borg52} and Marchenko \cite{Marchenko50}, and to modern
results by Simon \cite{BS1} and Gesztesy and Simon \cite{GS'00}.
The inverse problem on recovering the potential from the spectral
measure $d\rho$ was solved by Krein in \cite{Kr1,Kr2} and Gelfand
and Levitan in \cite{GL}.

\subsection{Quantum scattering inverse data.}

We consider the Schr\"odinger equation with the real-valued
potential $q\in L_{1+|x|}(\mathbb{R}_+)$
\begin{eqnarray}
&-\phi''+q(x)\phi=k^2 \phi,\quad x >0.\label{Schro_OP}
\end{eqnarray}
The solution $e(k,x)$ of the above equation is determined by the
condition
$$
\lim_{x\to\infty} e^{-ikx}e(k,x)=1.
$$
It admits the representation
\begin{equation*}
e(k,x)=e^{ikx}+\int_x^\infty K(x,t)e^{ikt}\,dt,
\end{equation*}
where the kernel $K(x,t)$ solves certain Goursat problem. The pair
$\{e(k,x),e(-k,x)\}$ forms a fundamental set of solutions when
$k\in \mathbb{R}$. Another solution to (\ref{Schro_OP})
$\varphi(k,x)$ is defined by the the conditions
\begin{equation*}
\varphi(k,0)=0,\quad \varphi_x(k,0)=1.
\end{equation*}
We set the notation $M(k)=e(0,k)$. Then $e$ and $\varphi$ when $k$
is on real axis are connected by
\begin{equation}
\label{Phi_e_conn}
-\frac{2ik\varphi(k,x)}{M(k)}=e(-k,x)-S(k)e(k,x),
\end{equation}
where the \emph{scattering matrix} is defined by
\begin{equation*}
S(k)=\frac{M(-k)}{M(k)}=\frac{1+\widehat K(0,-k)}{1+\widehat
K(0,k)},\quad k\in \mathbb{R}.
\end{equation*}
And on introducing the amplitude and phase of $M(k)$, we have:
\begin{eqnarray}
M(k)=A(k)e^{i\eta(k)},\quad A(k)=|M(k)|,\,\, \eta(k)=\arg{M(k)},\label{MK}\\
A(k)=A(-k),\quad \eta(k)=-\eta(-k).\label{MK1}
\end{eqnarray}

The problem (\ref{Schro_OP}) has a finite number of (negative)
eigenvalues $-k_1^2,\ldots, -k_n^2$, where $ik_l$ are zeroes of
the function $e(k,0)$, $l=1\ldots, n$. By $(C_j)^{-1}$ we denote
$(C_j)^{-1}=\int_0^\infty |e(ik_j,x)|^2\,dx$. Then the set of
functions
\begin{eqnarray*}
\{\varphi(k,x),\, k\in \mathbb{R}_+,\quad
\varphi_j(x)=e(ik_j,x),\, j=1,\ldots,n\}
\end{eqnarray*}
is a complete orthonormal set of eigenfunctions of the problem
(\ref{Schro_OP}). The Parseval identity has the form
\begin{equation*}
\delta(x-y)=\sum_{j=1}^n C_j^2\varphi_j(x)
\varphi_j(y)+\int_0^\infty
\varphi(x,k)\frac{1}{M(k)M(-k)}\varphi(y,k)k^2\,dk,
\end{equation*}
After we introduce notations (here $f\in L_2(\mathbb{R}_+)$)
\begin{eqnarray}
U(k):=\frac{1}{M(-k)M(k)}, \notag\\
\left(F^sf\right)(k)=\int_0^\infty f(x)\varphi(k,x)\,dx, \quad
\left(F^s_jf\right)=\int_0^\infty f(x)\varphi_j(x)\,dx,
\label{Fs_trans}
\end{eqnarray}
the Parseval equality for arbitrary $f,g\in L_2(\mathbb{R}_+)$
reads
\begin{equation*}
(f,g)_{L_2(\mathbb{R}_+)}=\sum_{j=1}^n C_j^2 \left(F^s_jf\right)
\left(F^s_jg\right)+\frac{2}{\pi}\int_0^\infty
\left(F^sf\right)(k)\left(F^sg\right)(k)U(k)k^2\,dk.
\end{equation*}
The set
\begin{equation*}
S_D:=\left\{(k_j,C_j)_{j=1}^n,\, S(k),\, k\in \mathbb{R}\right\}
\end{equation*}
is called the \emph{scattering data}. For the solution of the
inverse problem from $S_D$ see \cite{F63} and references therein.
It is important that the set $S_D$ determines the function $M(k)$
and thus, $U(k)$.

\subsection{Inverse acoustical scattering problem}

We consider the dynamical system associated with the (forward)
problem
\begin{equation}
\label{sc_eq} \left\{
\begin{array}l
u_{tt}-u_{xx}+qu=0,\quad  x>0,\,\,\, -\infty<t<x\\
u|_{t<-x}=0\\
\lim_{s\to\infty}u(s,\tau-s)=f(\tau), \quad \tau\geqslant 0,
\end{array}
\right.
\end{equation}
where $q\in C^\infty[0,\infty)$, $\operatorname{supp}q\subset
[0,a]$, $a<\infty$ is a \emph{potential}, $f$ is a \emph{control},
$u=u^f(x,t)$ is a solution (\emph{wave}).

Since $q|_{x>a}=0$, for large $x$'s the solution satisfies
$u_{tt}-u_{xx}=0$ and, hence, is a sum of two D'Alembert waves:
\begin{equation}
\label{Sc_sol_repr} u^f(x,t)|_{x>a}=f(x+t)+f^*(x-t),
\end{equation}
where the second summand describes the wave reflected by the
potential and outgoing to $x=\infty.$ Taking $f=\delta(t)$, one
can introduce a fundamental solution of the form
$u^\delta(x,t)=\delta(t+x)+w(x,t)$, which satisfies
\begin{equation}
\label{Sc_sol_repr1} u^\delta(x,t)|_{x>a}=\delta(x+t)+p(x-t)
\end{equation}
with $p\in C^\infty[0,\infty),$
$\operatorname{supp\,}p\subset[0,\infty]$. The Duhamel
representation $u^f=u^\delta*f$ holds for the classical solutions.
Note that $\operatorname{supp}f^*\subset [0,2a],$ so that the
reflected wave $f^*(x-t)$ in (\ref{Sc_sol_repr}) is compactly
supported on $t\leqslant x \leqslant \infty$ for any $t$.

An \emph{outer space} of the system (\ref{sc_eq})is the space of
controls $\mathcal{F}:=L_2(0,\infty)$. An \emph{inner space} is
$\mathcal{H}:=L_2(0,\infty)$ (of functions of $x$). A
\emph{control operator} $W:\mathcal{F}\to \mathcal{H}$ acts by the
rule
\begin{equation*}
(Wf)(x):=u^f(x,0),\qquad x\geqslant 0.
\end{equation*}
It maps $\mathcal{F}$ onto $\mathcal{H}$ isomorphically. These
facts are derived from the representation
\begin{equation}
\label{Sc_W_repr} (Wf)(x)=f(x)+\int_0^x w(x,-s)f(s)\,ds,\qquad
x\geqslant 0,
\end{equation}
A \emph{connecting operator} $C:\mathcal{F}\to\mathcal{F},$
\begin{equation*}
C:=W^*W
\end{equation*}
is a positive definite isomorphism in $\mathcal{F}$. It connects
the metrics of the outer and inner spaces:
\begin{equation}
\label{Sc_C}
(Cf,g)_{\mathcal{F}}=(Wf,Wg)_\mathcal{H}=\left(u^f(\,\cdot\,,0),u^g(\,\cdot\,,0)\right)_\mathcal{H}.
\end{equation}
\smallskip

\noindent A {\bf response operator} of the system (\ref{sc_eq}) is
$R:\mathcal{F}\to\mathcal{F},$
\begin{equation*}
(Rf)(\tau):=\lim_{s\to +\infty}u^f(s,s-\tau),\qquad \tau\geqslant
0.
\end{equation*}
For $f\in \mathcal{F}$ vanishing at $\infty$, by
(\ref{Sc_sol_repr}), this limit is $f^*(\tau).$ Hence we get
\begin{equation}
\label{Sc_R_repr} (Rf)(\tau)=\int_0^\infty
p(\tau+s)f(s)\,ds,\qquad \tau\geqslant 0.
\end{equation}
Here $p$, the \emph{acoustical response function} could be
determined as a response to delta function:
\begin{equation*}
p(\tau)=\lim_{s\to +\infty}u^\delta(s,s-\tau),\quad \tau>0.
\end{equation*}
The inverse acoustical scattering problem is to recover potential
$q|_{x\geqslant 0}$ by given response operator $R$ (or what is
equivalent, from acoustical response $p|_{t\geqslant 0}$).
\begin{theorem}
\label{Sc_teor_C_repr} The equality
\begin{equation}
\label{Sc_C_repr} C=\mathbb{I}+R
\end{equation}
holds, where $\mathbb{I}$ is the identity operator in
$\mathcal{F}$.
\end{theorem}

A natural setup of a \emph{control problem} for the system
(\ref{sc_eq}) is by given $y\in \mathcal{H}$ to find $f\in
\mathcal{F}$ such that $u^f(\,\cdot\,,0)=y$. This problem is
equivalent to the equation $Wf=y$, which has a unique solution
$f=W^{-1}y\in \mathcal{F}$ due to (\ref{Sc_W_repr})).

Let us consider a \emph{special control problem}: take $y$, which
satisfies
\begin{equation}
\label{Sct_y}
\left\{
\begin{array}l
 -y''+qy=k^2y,\qquad x>0,\quad \\
y|_{x>a}=e^{ikx}.
\end{array}
\right.
\end{equation}
\begin{theorem}
\label{Sc_control_teor} The control $f=W^{-1}y,$ which solves the
special CP, is the solution of the equation $Cf=e^{ik(\cdot)}$ in
$\mathcal{F}$, i.e., it satisfies
\begin{equation}
\label{Sc_C_eq} f(\tau)+\int_0^\infty
r(\tau+s)f(s)\,ds=e^{ik\tau},\qquad \tau\geqslant 0.
\end{equation}
\end{theorem}

Writing (\ref{Sc_C_eq}) in the form $W^*Wf=e^{ik(\cdot)},$ with
regard to $Wf=y,$ we have
\begin{equation}
\label{AcSc_transf} W^*y(\cdot,k)=e^{ik(\cdot)}.
\end{equation}
Hence, $W^*$ is a {\it transformation operator}, which maps the
solution $y(\cdot,k)$ of (\ref{Sct_y}) to the solution
$e^{ik\cdot}$ of the unperturbed problem.

\section{Kernel of the connecting operator $C^T$}

\subsection{The spectral function of Levitan and the kernel of the connecting operator.}

Here we derive the spectral representation of the connecting
operator (\ref{Ct}), (\ref{r-c2}) following \cite{AM}.

We take a Fourier transform (\ref{Fourier_int_2}) of
$u^f(\cdot,T)$ and use the transformation property
(\ref{WT_wv_tr}) of $\left(W^T\right)^*$ :
\begin{eqnarray}
\label{DB_Sch}\left(Fu^f(\cdot,T)\right)(\mu)=\int_{-\infty}^\infty
u^f(x,T)\varphi(x,\mu)\,dx=(W^Tf,\varphi(\cdot,\mu))_{\mathcal{H}^T}\\
=(f,\left(W^T\right)^*\varphi(\cdot,\mu))_{\mathcal{H}^T}=\int_0^T\frac{\sin{\sqrt{\mu}s}}{\sqrt{\mu}}f(T-s)\,ds.\notag
\end{eqnarray}
Let  $f,g\in \mathcal{F}^T$. Using (\ref{Fourier_int_1}) and
(\ref{Fourier_int_2}), we rewrite
$\left(C^Tf,g\right)_\mathcal{F^T}$ as
\begin{equation}
\label{C_T_1} \left(C^Tf,g\right)_{\mathcal{F}^T}=\int_0^T
u^f(x,T)u^g(x,T)\,dx=\int_{-\infty}^{\infty}
(Fu^f)(\lambda,T)(Fu^g)(\lambda,T)\,d\rho(\lambda).
\end{equation}
Making use of (\ref{DB_Sch}), we can rewrite (\ref{C_T_1}) as
\begin{equation}
\label{CT_1}
\left(C^Tf,g\right)_{\mathcal{F}^T}=\int_{-\infty}^\infty\,\int_0^T\,\int_0^T\,
\,\frac{\sin{\sqrt{\lambda}(T-t)}\sin{\sqrt{\lambda}(T-s)}}{\lambda}f(t)g(s)\,dt\,ds\,d\rho(\lambda).
\end{equation}
Now we make use of the $\sin$ transform: for all $h, j\in
L^2(\mathbb{R}_+)$
\begin{eqnarray*}
\widehat h(\lambda)=\int_0^\infty
h(x)\frac{\sin{(\sqrt{\lambda}x)}}{\sqrt{\lambda}}\,dx,\quad
h(x)=\int_0^\infty \widehat h(\lambda)\sin{(\sqrt{\lambda}x)}\,d\left(\frac{2}{3\pi}{\lambda}^\frac{3}{2}\right),\\
\int_0^\infty h(x)j(x)\,dx=\int_0^\infty \widehat
h(\lambda)\widehat
j(\lambda)\,d\left(\frac{2}{3\pi}{\lambda}^\frac{3}{2}\right).
\end{eqnarray*}
Let us extend the functions $f$ and $g$ to the whole real axis
 setting $f(t)=g(t)=0$ for $t>T$ and $t<0$ and use the
notation $f_T(s)=f(T-s)$. Then we can rewrite
\begin{eqnarray}
\int_0^T f(t)g(t)\,dt=\int_0^\infty f(T-s)g(T-s)\,ds =
\int_0^\infty \widehat f_T(\lambda)\widehat
g_T(\lambda)\,d\left(\frac{2}{3\pi}{\lambda}^\frac{3}{2}\right)\notag\\
=\int_0^\infty\,\int_0^T\, \int_0^T\,
\frac{\sin{\sqrt{\lambda}(T-t)}\sin{\sqrt{\lambda}(T-s)}}{\lambda}f(t)g(s)\,dt\,ds\,d\left(\frac{2}{3\pi}{\lambda}^\frac{3}{2}\right).\label{FG}
\end{eqnarray}
On introducing the function
\begin{eqnarray*}
\sigma(\lambda)=\left\{\begin{array}l
\rho(\lambda)-\frac{2}{3\pi}{\lambda}^{\frac{3}{2}},\quad \lambda\geqslant 0, \\
\rho(\lambda),\quad \lambda <0,
\end{array}
\right.
\end{eqnarray*}
we can rewrite (\ref{CT_1}) using (\ref{FG}) and counting that for
fixed $n$ we we can change the order of integration:
\begin{eqnarray}
\label{C_T_5}
\lim_{n\to\infty}\int_0^T\,\int_0^T\,\int_{-\infty}^n
\,\frac{\sin{\sqrt{\lambda}(T-t)}\sin{\sqrt{\lambda}(T-s)}}{\lambda}\,d\sigma(\lambda)f(t)g(s)\,dt\,ds \\
= \int_0^T\,\int_0^T\,c^T(s,t)f(t)g(s)\,dt\,ds.\notag
\end{eqnarray}
Let us introduce the function
\begin{equation}
\label{Psi_n} \Psi_n(t,s):=\int_{-\infty}^n
\,\frac{\sin{\sqrt{\lambda}(T-t)}\sin{\sqrt{\lambda}(T-s)}}{\lambda}\,d\sigma(\lambda)
\end{equation}
Since $f,g$ are arbitrary functions from $\mathcal{F}^T$, we can
deduce from (\ref{C_T_5}) that
\begin{equation*}
\Psi_n(t,s)\longrightarrow_{n\to\infty} c^T(t,s),\quad
\text{weekly in $L_2\left((0,T)^2\right)$}.
\end{equation*}
To strengthen the result on the convergence we need the theorem of
Levitan \cite{L3} on the convergence of spectral functions:
\begin{theorem}
The sequence of functions
\begin{equation} \label{Levitan_sf}
\Phi_n(t,s)=\int_{-\infty}^n
\f(t,\lambda)\f(s,\lambda)\,d\rho(\lambda)-\int_0^n
\frac{\sin{\sqrt{\lambda}t}\sin{\sqrt{\lambda}s}}{\lambda}\,d\left(\frac{2}{3\pi}\lambda^{\frac{3}{2}}\right),
\end{equation}
converges uniformly on every bounded set to a differentiable
outside the diagonal function, as $n$ tends to infinity.
\end{theorem}
Applying operators $({\bf I}_s+{\bf L}_s)({\bf I}_t+{\bf L}_t)$
(see (\ref{Inv_trans_wave})) to (\ref{Levitan_sf}) we have:
\begin{eqnarray}
\label{L_eq}
({\bf I}_s+{\bf L}_s)({\bf I}_t+{\bf L}_t) \Phi_n(t,s)=\Psi_n(s,t) \\
-
\int_0^n\left(\int_0^tL(t,\tau)\frac{\sin{\sqrt{\lambda}\tau}}{\sqrt{\lambda}}\,d\tau\right)\frac{\sin{\sqrt{\lambda}s}}{\sqrt{\lambda}}\,
d\left(\frac{2}{3\pi}\lambda^{\frac{3}{2}}\right) \notag\\
-
\int_0^n\left(\int_0^sL(s,\tau)\frac{\sin{\sqrt{\lambda}\tau}}{\sqrt{\lambda}}\,d\tau\right)\frac{\sin{\sqrt{\lambda}t}}{\sqrt{\lambda}}\,
d\left(\frac{2}{3\pi}\lambda^{\frac{3}{2}}\right) \notag\\
-
\int_0^n\left(\int_0^tL(t,\tau)\frac{\sin{\sqrt{\lambda}\tau}}{\sqrt{\lambda}}\,d\tau\right)
\left(\int_0^sL(s,\tau)\frac{\sin{\sqrt{\lambda}\tau}}{\sqrt{\lambda}}\,d\tau\right)\,d\left(\frac{2}{3\pi}\lambda^{\frac{3}{2}}\right).\notag
\end{eqnarray}
The sum of the last three terms in the right hand side of the
above expression converges to
$-L(s,t)-L(t,s)-\int_0^{\min{\{s,t\}}}L(s,\tau)L(t,\tau)\,d\tau$.
This fact and the convergence of the left hand side of
(\ref{L_eq}) imply that
\begin{equation}
\label{CT_spectr_resp} \Psi_n(t,s)\longrightarrow_{n\to\infty}
c^T(t,s)=\int_{-\infty}^\infty
\,\frac{\sin{\sqrt{\lambda}(T-t)}\sin{\sqrt{\lambda}(T-s)}}{\lambda}\,d\sigma(\lambda),
\end{equation}
uniformly on every compact set in $\mathbb{R}^2$.

The estimates on regularized spectral function $\Phi_n(t,s)$
receive a lot of attention, to mention \cite{Sad} and literature
cited therein. We believe that the connection of the regularized
spectral function $\Phi_n$ with the kernel of the connecting
operator $C^T$ allows one to extends some of the results to
different dynamical systems, for example to vector Schr\"odinger
system, Dirac system, canonical systems.

\section{ On the connection of the spectral, dynamical and scattering data.}

\subsection{Weyl function and response function}

We now demonstrate the connection between the response function
$r\left( s\right) $ and the Titchmarsh-Weyl $m$-function. A
connection between spectral and time-domain data is widely used in
inverse problems, see, e.g., \cite{B01JII,B03,KKLM} where the
equivalence of several types of boundary inverse problems is
discussed.

Let $f \in C_0^{\infty}(0,\infty)$ and
\begin{equation*}
\widehat{f}(k):=\int_0^{\infty} f(t)\,e^{-kt}\,dt\,
\end{equation*}
be its Laplace transform. Function $\widehat{f}(k)$ is well
defined for $k \in \mathbb{C}$. Going in (\ref{wave_eqnGL}) over
to the Laplace transforms, one has
\begin{equation*}
\left\{ \begin{array}l
-\widehat{u}_{xx}(x,k)+q(x)\widehat{u}(x,k)
=-k^{2}\widehat{u}(x,k),\\
\widehat{u}(0,k) =\widehat{f}(k),
\end{array}
\right.
\end{equation*}
and
\begin{equation*}
\widehat{(Rf)}(k)=\widehat{u}_{x}(0,k).
\end{equation*}
respectively. The values of the function $\widehat{u} (0,k) $ and
its first derivative at the origin, $\widehat{u}_x (0,k), $ are
related through the Titchmarsh-Weyl m-function
\begin{equation*}
\widehat{u}_x (0,k)=m(-k^2)\widehat{f}(k)\,.
\end{equation*}
Therefore,
\begin{equation}
\label{2.7} \widehat{(Rf)}(k)=m(-k^2)\widehat{f}(k)\,,
\end{equation}
and thus the spectral and dynamic Dirichlet-to-Neumann maps are in
one-to-one correspondence. Taking the Laplace transform of
(\ref{react_rep}) we get
\begin{equation}
\widehat{(Rf)}(k) =-k\widehat{f}(k)+\widehat{r}(k)\widehat{f}(k).
\label{2.8}
\end{equation}%
In \cite{AMR} the authors have shown that, if the potential
belongs to the class
$$
q\in
l^{\infty }\left( L^{1}\left( \mathbb{R}_{+}\right) \right)
:=\left\{ q:\int_{n}^{n+1}\left\vert q\left( x\right) \right\vert
\,dx\in l^{\infty }\right\} .
$$
with the norm defined by  $\left\vert \left\vert q\right\vert
\right\vert =\sup \int_{x}^{x+1}\left\vert q\left( s\right)
\right\vert ds<\infty$ . Then (\ref{2.7}) and (\ref{2.8}) imply
\begin{equation}
m(-k^{2})=-k+\int_{0}^{\infty }e^{-k\alpha }r(\alpha )\,d\alpha \,
, \label{2.5}
\end{equation}
with the integral in  (\ref{2.5}) is absolute convergent for $z=-k^{2}$ where $%
\operatorname{Re}k>2\max \{\sqrt{2\left| \left| q\right| \right|
},e\left| \left| q\right| \right| \}.$

We notice that it was shown in \cite{BS1} that there exists a
unique real valued function $A \in L^{1}_{loc}\left(
\mathbb{R}_{+}\right) $ (the\emph{\ }$A-$\emph{amplitude}) such
that
\begin{equation}
\label{A_amp} m(-k^{2})=-k-\int_{0}^{\infty }A(t )e^{-2t k}\,dt\,.
\end{equation}
The absolute convergence of the integral was proved for $q\in
L^{1}\left( \mathbb{R}_{+}\right) $ and $q\in L^{\infty }\left( \mathbb{R}%
_{+}\right) $ in \cite{GS'00} for sufficiently large $\Re{k}$.
Clearly, the $A-$amplitude and response function are connected by
\begin{equation*}
A(t)=-2r(2t).
\end{equation*}

\subsection{Response function and spectral measure}

Using the representation for $c^T(t,s)$ (\ref{CT_spectr_resp}), we
can derive the formula for the response function: \begin{theorem}
\label{Th_Resp} The representation for the response function $r$
\begin{equation}
\label{Resp_mes_conn} r(t)=\int_{-\infty}^\infty
\frac{\sin{\sqrt{\lambda}t}}{\sqrt{\lambda}}\,d\sigma(\lambda),\,
\end{equation}
holds for almost all $t \in [0,+\infty)$.
\end{theorem}
\begin{proof}
Let us note that according to Theorem
\begin{equation}
\label{Ph_i} \Phi(s,t)=\int_{-\infty}^\infty
\frac{\sin{\sqrt{\lambda}t}\sin{\sqrt{\lambda}s}}{\lambda}\,d\sigma(\lambda)=c^T(T-t,T-s),\quad
t,s\in [0,T].
\end{equation}
Using (\ref{c_t}), we have
\begin{equation}
\label{Phi_1}
c^T(T-t,T-s)=\frac{1}{2}\int_{|t-s|}^{t+s}r(\tau)\,d\tau,\quad
t,s\in [0,T].
\end{equation}
The integral in (\ref{Ph_i}) can be rewritten as
\begin{eqnarray}
\Phi(s,t)=\frac{1}{2}\int_{-\infty}^\infty
\frac{(\cos{\sqrt{\lambda}(s+t)}-1-(\cos{\sqrt{\lambda}|s-t|}-1)}{\lambda}\,d\sigma(\lambda)=\label{Phi_2}\\
\frac{1}{2}\int_{-\infty}^\infty\int_{|t-s|}^{t+s}\frac{\sin{\sqrt{\lambda}\theta}}{\sqrt{\lambda}}\,d\theta\,d\sigma(\lambda),\quad
t,s\in [0,T].\notag
\end{eqnarray}
Equating the expressions in  (\ref{Phi_1}) and (\ref{Phi_2}) for
$t=s$ we get
\begin{equation}
2c^T(T-t,T-t)=\int_0^{2t}r(\tau)\,d\tau=
\int_{-\infty}^\infty\int_0^{2t}\frac{\sin{\sqrt{\lambda}\theta}}{\sqrt{\lambda}}\,d\theta\,d\sigma(\lambda),\quad
t\in [0,T].
\end{equation}
According to (\ref{react_rep}),  $r\in L^1(0,T)$, so we can use
the Lebesgue theorem and differentiate the last equation. We
obtain the following equality almost everywhere on $[0,2T]$:
\begin{equation*}
r(t)=\int_{-\infty}^\infty\frac{\sin{\sqrt{\lambda}t}}{\sqrt{\lambda}}\,d\sigma(\lambda).
\end{equation*}
Since the parameter $T$  can be chosen arbitrary large, the last
formula proves the statement of the proposition.
\end{proof}
A different proof is given in \cite{AM}, see also \cite{BS1},
where the regularized version of (\ref{Resp_mes_conn}) was
derived.

The finite speed of the wave propagation in (\ref{wave_eqnGL})
implies the local nature of the response function $r(t)$: the
values of $r(t)$, $t\in [0,2T]$ are determined by the potential
$q(x)$, $x\in [0,T]$. That is why if we are interested in the
representation of $r(t)$ on the interval $t\in [0,2T]$, we can
replace  in (\ref{Resp_mes_conn}) the regularized spectral
function $\sigma(\lambda)$ by any of the following functions:
\begin{equation*}
\sigma_{\operatorname{tr}}(\lambda)=\left\{\begin{array}l
\rho_{\operatorname{tr}}(\lambda)-\frac{2}{3\pi}{\lambda}^{\frac{3}{2}},
\quad \lambda\geqslant 0,\\
\rho_{\operatorname{tr}}(\lambda),\quad \lambda
<0,\end{array}\right.,\quad
\sigma_{d}(\lambda)=\left\{\begin{array}l
\rho_{d}(\lambda)-\rho_0(\lambda),\quad \lambda\geqslant 0,\\
\rho_{d}(\lambda),\quad \lambda <0.\end{array}\right.
\end{equation*}
Here $\rho_{tr}$ is the spectral function corresponding to the
truncated potential: $q_T(x)=q(x)$ for $0 \leq x \leq T$ and
$q_T(x)=\tilde q(x)$ for $x>T$ with arbitrary locally integrable
$\tilde q$; $\rho_d(\lambda)$ is the spectral function associated
to the discrete problem on the interval $(0,T)$ with the potential
$q_d(x)=q(x)$, $x\in [0,T]$ and $\rho_0(\lambda)$ is the spectral
function associated to the discrete problem on $[0,T]$ with zero
potential. (Any self-adjoint boundary condition can be prescribed
at $x=T.$)

We notice that the function $\Phi(t)=\int_0^tr(\tau)\,d\tau$, in
accordance with (\ref{Resp_mes_conn}) is given by
\begin{equation*}
\Phi(t)=\int_{-\infty}^\infty
\frac{1-\cos{\sqrt{\lambda}\tau}}{\lambda}\,d\sigma(\lambda),\quad
0<t<2T.
\end{equation*}
has been used by Krein in \cite{Kr1,Kr2} as an inverse data.

\subsection{Quantum scattering data and response function}

Using the representation for $c^T(t,s)$ obtained in
(\ref{CT_spectr_resp}), we can derive the formula for the response
function:
\begin{theorem}
\label{Th_Resp1} The representation for the response function $r$
in terms of scattering data:
\begin{equation}
\label{Resp_mes_conn1} r(t)=\sum_{j=1}^n
C_j^2\frac{\sin{k_jt}}{k_j}+\frac{2}{\pi}\int_0^\infty \sin{k
t}\left(U(k)-1\right)k\,dk
\end{equation}
holds for almost all $t \in [0,+\infty)$.
\end{theorem}
\begin{proof}
Let us take arbitrary $f, g\in \mathcal{F}^T$ and consider the
connecting operator $C^T$ (\ref{Ct})
\begin{equation}
\label{C_T}
\left(C^Tf,g\right)_{\mathcal{F}^T}=\left(u^f(\cdot,T),u^g(\cdot,T)\right)_{\mathcal{H}^T}.
\end{equation}
Where $u^f$ is solutions to the wave equation (\ref{wave_eqnGL})
with the control $f$. Rewriting (\ref{C_T}) using the Parseval
identity, we obtain
\begin{eqnarray*}
\left(C^Tf,g\right)_{\mathcal{F}^T}=\sum_{j=1}^n C_j^2
F^s_j(u^f(\cdot,T))
F^s_j(u^g(\cdot,T))\\
+\frac{2}{\pi}\int_0^\infty
\left(F^su^f(\cdot,T)\right)(k)\left(F^su^g(\cdot,T)\right)(k)U(k)k^2\,dk
\end{eqnarray*}
Using the transformation property (\ref{WT_wv_tr}) of
$\left(W^T\right)^*$ yields
\begin{eqnarray*}
\left(F^su^f(\cdot,T)\right)(k)= \int_0^T
\varphi(x,k)u^f(x,T)\,dx=
\left(\varphi(x,k),W^Tf\right)_{H^T}\\
=\left(\left(W^T\right)^*\varphi(x,k),f\right)_{H^T}=\int_0^T
f(T-s)\frac{\sin{ks}}{k}\,ds.
\end{eqnarray*}
Similarly,
\begin{eqnarray*}
\left(F^s_ju^f(\cdot,T)\right)(k)=\int_0^T
f(T-s)\frac{\sin{k_js}}{k_j}\,ds.
\end{eqnarray*}
Using these observations we get an equivalent expression for
(\ref{C_T}):
\begin{eqnarray}
\left(C^Tf,g\right)_{\mathcal{F}^T}=\sum_{j=1}^n C_j^2
\int_0^T\int_0^T
\frac{\sin{k_j(T-s)}}{k_j}\frac{\sin{k_j(T-\tau)}}{k_j}f(s)g(\tau)\,ds\,d\tau\notag\\
+\frac{2}{\pi}\int_0^\infty\int_0^T\int_0^T
\frac{\sin{k(T-s)}}{k}\frac{\sin{k(T-\tau)}}{k}
U(k)k^2f(s)g(\tau)\,ds\,d\tau\,dk \label{C_sc1}
\end{eqnarray}
Using the representation for $C^T$ (\ref{r-c2}), (\ref{c_t})
and
\begin{equation*}
\int_0^Tf(t)g(t)\,dt=\frac{2}{\pi}\int_0^\infty\int_0^T\int_0^T
\frac{\sin{k(T-s)}}{k}\frac{\sin{k(T-\tau)}}{k}
k^2f(s)g(\tau)\,ds\,d\tau\,dk,
\end{equation*}
we can rewrite (\ref{C_sc1}) as
\begin{eqnarray}
\label{C_sc3} \int_0^T\int_0^T
c^T(t,s)f(t)g(s)\,dt\,ds\\=\int_0^T\int_0^T \sum_{j=1}^n
C_j^2\frac{\sin{k_j(T-s)}}{k_j}\frac{\sin{k_j(T-\tau)}}{k_j}f(s)g(\tau)\,ds\,d\tau\notag\\
+\frac{2}{\pi}\int_0^\infty\int_0^T\int_0^T
\frac{\sin{k(T-s)}}{k}\frac{\sin{k(T-\tau)}}{k}
(U(k)-1)k^2f(s)g(\tau)\,ds\,d\tau\,dk.\notag
\end{eqnarray}
We notice that it is possible to change the order of integration
in the last integral in (\ref{C_sc3}) due to results on
convergence from \cite{F63}. The latter observation leads to the
representation for the kernel $c(t,x)$:
\begin{eqnarray*}
c(t,x)=\sum_{j=1}^n C_j^2
\frac{\sin{k_j(T-x)}}{k_j}\frac{\sin{k_j(T-t)}}{k_j}+\\
\frac{2}{\pi}\int_0^\infty
\frac{\sin{k(T-t)}}{k}\frac{\sin{k(T-x)}}{k}
\left(U(k)-1\right)k^2\,dk\notag
\end{eqnarray*}
We have that on the one hand (\ref{Phi_1}), and on the other hand
\begin{eqnarray}
\label{e_2} c(T-t,T-x)=\sum_{j=1}^n C_j^2
\frac{\sin{k_jx}}{k_j}\frac{\sin{k_jt}}{k_j}+\\
\frac{2}{\pi}\int_0^\infty \frac{\sin{k t}}{k}\frac{\sin{kx}}{k}
\left(U(k)-1\right)k^2\,dk\notag
\end{eqnarray}
Note that
\begin{equation}
\label{e_3} \frac{\sin{k
t}}{k}\frac{\sin{kx}}{k}=\frac{1}{2}\int_{|t-x|}^{t+x}\frac{\sin{k\theta}}{k}\,d\theta,
\end{equation}
Using (\ref{Phi_1}), (\ref{e_2}) and (\ref{e_3}) we arrive at
(\ref{Resp_mes_conn1}). This formula but in weaker form was
derived in \cite{GS'00}.
\end{proof}

\subsection{Scattering matrix and acoustical response function.}

Let $f \in C_0^{\infty}(0,\infty)$ and
\begin{equation*}
\left(Ff\right)(k):=\int_{-\infty}^{\infty} f(t)\,e^{ikt}\,dt\,
\end{equation*}
be its Fourier transform. Function $\widehat{f}(k)$ is well
defined for $k \in \mathbb{C}$. Going in (\ref{sc_eq}) with
control $f=\delta$ over to the Fourier transforms, one has
\begin{equation*}
\left\{
\begin{array}l
-\widehat{u}^\delta_{xx}(x,k)+q(x)\widehat{u}^\delta(x,k)
=k^{2}\widehat{u}^\delta(x,k), \\
\widehat{u}(0,k)=0,
\end{array}
\right.
\end{equation*}
On the other hand, for $x>a$ one has the representation
(\ref{Sc_sol_repr1}) for $u^\delta$, applying Fourier transform to
it, we get
\begin{equation}
\label{Acus_Four}
\widehat{u}^\delta(x,k)=e^{-ikx}+\int_{-\infty}^\infty
p(s)e^{-iks}\,ds\, e^{ikx},\quad x>a.
\end{equation}
Comparing (\ref{Acus_Four}) with (\ref{Phi_e_conn}), we conclude
that scattering matrix and acoustical response are connected by
\begin{equation}
\label{Scatter_matr}
S(k)=-\int_{-\infty}^\infty
p(s)e^{-iks}\,ds=-2\pi \left(F^{-1}p\right)(k).
\end{equation}
We study the connecting operator for the acoustical problem
(\ref{Sc_C}). Using the transformation $F^s$ (\ref{Fs_trans}), we
rewrite (\ref{Sc_C}) as
\begin{eqnarray}
(Cf,g)_{\mathcal{F}}=(Wf,Wg)_\mathcal{H}=\sum_{j=1}^n C_j^2
\left(F^s_jWf\right)
\left(F^s_jWg\right)\label{CT_sc_dev}
\\+\frac{2}{\pi}\int_0^\infty
\left(F^sWf\right)(k)\left(F^sWg\right)(k)U(k)k^2\,dk\notag
\end{eqnarray}
Let us evaluate using (\ref{Phi_e_conn}) and transformation
property of $W^*$ (\ref{AcSc_transf}):
\begin{eqnarray}
\left(F^sWf\right)(k)=\int_0^\infty
\left(Wf\right)(x)\varphi(x,k)\,dx\label{FS}\\
=\int_0^\infty
\left(Wf\right)(x)\frac{i}{2k}\left(M(k)e(x,-k)-M(-k)e(x,k)\right)\,dx\notag\\
=\frac{i}{2k}\left(f,M(k)e^{-ik\cdot}-M(-k)e^{ik\cdot}\right)\notag
\end{eqnarray}
\begin{eqnarray}
\left(F^s_jWf\right)(k)=\int_0^\infty
\left(Wf\right)(x)\varphi_j(x)\,dx=\left(Wf,e(ik_j,\cdot)\right)\label{FS1}\\
=\left(f,W^*e(ik_j,\cdot)\right)=\left(f,e^{-k_j\cdot}\right)=:f_j.\notag
\end{eqnarray}
We continue evaluate (\ref{CT_sc_dev}) using (\ref{FS}),
(\ref{FS}) and (\ref{MK}), (\ref{MK1}):
\begin{eqnarray}
(Cf,g)_{\mathcal{F}}=\sum_{j=1}^n C_j^2
f_jg_j\notag
\\+\frac{2}{\pi}\int_0^\infty
\left(f,\frac{e^{i(\eta(k)-k\cdot)}-e^{i(k\cdot-\eta(k))}}{2i}\right)\left(g,\frac{e^{i(\eta(k)-k\cdot)}-e^{i(k\cdot-\eta(k))}}{2i}\right)\,dk\notag\\
=\sum_{j=1}^n C_j^2 f_jg_j+\frac{2}{\pi}\int_0^\infty
\left(f,\sin{(k\cdot-\eta(k))}\right)\left(g,\sin{(k\cdot-\eta(k))}\right)\,dk
\end{eqnarray}
In the dynamical representation (\ref{Sc_C}), (\ref{Sc_C_eq}):
\begin{equation*}
(Cf,g)_{\mathcal{F}}=\frac{2}{\pi}\int_0^\infty
\left(f,\sin{(k\cdot)}\right)\left(g,\sin{(k\cdot)}\right)\,dk+\int_0^\infty
\int_0^\infty p(t+s)f(t)g(s)\,dt\,ds.
\end{equation*}
Then we can write:
\begin{eqnarray*}
\int_0^\infty \int_0^\infty p(t+s)f(t)g(s)\,dt\,ds=
\int_0^\infty\int_0^\infty \sum_{j=1}^n
C_j^2e^{-k_j(x+y)}f(x)g(y)\,dx\,dy\\
+\int_0^\infty\int_0^\infty f(x)g(y) \frac{2}{\pi}\int_0^\infty
\sin{(kx-\eta(k))}\sin{(ky-\eta(k))}\,dk\\
-\int_0^\infty\int_0^\infty f(x)g(y) \frac{2}{\pi}\int_0^\infty
\sin{(kx)}\sin{(ky)}\,dk
\end{eqnarray*}
by the trigonometry,
\begin{eqnarray*}
\int_0^\infty \int_0^\infty p(t+s)f(t)g(s)\,dt\,ds=
\int_0^\infty\int_0^\infty \sum_{j=1}^n
C_j^2e^{-k_j(x+y)}f(x)g(y)\,dx\,dy\\
+\int_0^\infty\int_0^\infty f(x)g(y) \frac{1}{\pi}\int_0^\infty
\left(\cos{k(x+y)}-\cos{(k(x+y)-2\eta(k))}\right)\,dk
\end{eqnarray*}
from where we deduce the representation for $p$:
\begin{eqnarray}
p(t)=\sum_{j=1}^n
C_je^{-k_jt}+\frac{1}{\pi}\int_0^\infty\left(\cos{kt}-\cos{(kt-2\eta(k))}\right)\,dk,\label{ac_resp1}\\
p(t)=\sum_{j=1}^n
C_je^{-k_jt}+\frac{2}{\pi}\int_0^\infty\left(\sin{\eta(k)}\sin{(\eta(k)-kt)}\right)\,dk,\label{ac_resp2}
\end{eqnarray}
where the right hand sides are understood as generalized
functions.

Notice that the acoustical response in (\ref{Scatter_matr}) plays
the same role for the scattering matrix as response function (or
$A-$amplitude) plays for Weyl function in (\ref{2.5}) and
(\ref{A_amp}).

We think that the result of the convergence of integrals in
(\ref{ac_resp1}), (\ref{ac_resp2}) can be improved, we are planing
to come return to this question elsewhere in the framework of
studying the inverse acoustical scattering problem for the system
(\ref{sc_eq}) with a potential from a wider class: nonsmooth, and
sufficiently fast decreasing at infinity.

\noindent{\bf Acknowledgments}

The research of Victor Mikhaylov was supported in part by NIR
SPbGU 11.38.263.2014 and RFBR 14-01-00535. Alexandr Mikhaylov was
supported by RFBR 14-01-00306; A. S. Mikhaylov and V. S. Mikhaylov
were partly supported by VW Foundation program "Modeling,
Analysis, and Approximation Theory toward application in
tomography and inverse problems."

\end{document}